\documentclass{article}

\usepackage{delarray,verbatim,enumerate,a4wide}
\usepackage{amsmath,amsthm,amstext,amsbsy,amssymb,amsfonts,amscd}
\usepackage[all]{xy}
\usepackage[frenchb]{babel}
\usepackage[T1]{fontenc}
\usepackage[utf8]{inputenc}
\usepackage[small]{titlesec}
\usepackage{scalerel}
\usepackage{color}
\definecolor{vert}{rgb}{0.1,0.5,0.2}
\usepackage[colorlinks,linkcolor=blue,urlcolor=blue,citecolor=vert,linktocpage]{hyperref}

\usepackage{calligra}
\DeclareFontShape{T1}{calligra}{m}{n}{<->s*[0.95]callig15}{}
\DeclareMathAlphabet{\mathscr}{T1}{calligra}{m}{n}

\newtheorem{Th}{Théorème}
\newtheorem{Lem}[Th]{Lemme}
\newtheorem{Prop}[Th]{Proposition}
\newtheorem{Cor}[Th]{Corollaire}

\newtheorem{Sco}[Th]{Scolie}
\newtheorem{Def} [Th]{Définition}

\newtheorem*{ThA}{Théorème A}
\newtheorem*{ThB}{Théorème B}

\newtheorem*{Conv}{Convention}

\def\Preuve{\noindent {\it Preuve.~}}
\def\PreuveTh{\noindent {\it Preuve du Théorème.~}}
\def\Remarque{\noindent {\it Remarque.~}}

		\def\QQ{\mathbb Q}	
\def\NN{\mathbb N}	\def\ZZ{\mathbb Z}	
\def\F2{\mathbb{F}_2}	\def\Z2{\mathbb{Z}_2}		
\def\Zl{\mathbb{Z}_\ell} 	\def\Ql{\mathbb{Q}_\ell}	\def\Fl{\mathbb{F}_\ell}

 				\def\U{\mathcal  U}	\def\F{\mathcal  F}
\def\J{\mathcal  J}  		\def\R{\mathcal  R}	
\def\Dl{\mathcal  D\ell} 	\def\Pl{\mathcal  P\si{\!}\ell}  	\def\Cl{\mathcal  C \!\ell}
\def\E{\mathcal  E}				\def\	

\def\q{{\mathfrak q}}	\def\p{{\mathfrak p}}		\def\a{{\mathfrak a}} \def\b{{\mathfrak b}}
		\def\d{{\mathfrak d}}	\def\m{{\mathfrak m}}

\def\wi{\widetilde}		
	
	\def\deg{\operatorname{deg}}
\def\Gal{\operatorname{Gal}}		\def\Rad{\operatorname{Rad}}
\def\Ker{\operatorname{Ker}}		

\newcommand\scale[2]{\vstretch{#1}{\hstretch{#1}{#2}}}

\newcommand\si[1]{\scale{.6}{#1}}
\newcommand\ph{{\phantom{*}}}
\newcommand\ab{{\scale{.8}{ab}}}

\makeatletter
\newcommand*\wt[2][0.2ex]{%
        \begingroup
        \mathchoice{\wt@helper{#1}{#2}{\displaystyle}{\textfont}}
                   {\wt@helper{#1}{#2}{\textstyle}{\textfont}}
                   {\wt@helper{#1}{#2}{\scriptstyle}{\scriptfont}}
                   {\wt@helper{#1}{#2}{\scriptscriptstyle}{\scriptscriptfont}}%
        \endgroup
        #2%
}
\newcommand*\wt@helper[4]{%
        \def\currentfont{\the#41}%
        \def\currentskewchar{\char\the\skewchar\currentfont}%
        \setbox\tw@\hbox{\currentfont$#2$\currentskewchar}%
        \dimen@ii\wd\tw@
        \setbox\tw@\hbox{\currentfont$#2${}\currentskewchar}%
        \advance\dimen@ii-\wd\tw@
        \rlap{\raisebox{-#1}{$\m@th#3\kern\dimen@ii\widetilde{\phantom{#2}}$}}%
}
\makeatother

\def\wE{\,\wt[0.2ex]{\!\mathcal E}}		\def\wU{\wt[0.2ex]{\mathcal U}}
\def\wJ{\,\wt[0.2ex]{\!\mathcal J}}	\def\wCl{\wt[0.1ex]{\mathcal C\!\ell}} \def\wDl{\wt[0.2ex]{\mathcal D_{\si{\!}}\ell}}
\def\wT{\,\wt[0.2ex]{\!\mathcal T}}		\def\wH{\,\wt[0.2ex]{H}}			\def\wK{\,\wt[0.2ex]{K}}

\def\p{\mathfrak p}		\def\q{\mathfrak q}			\def\r{\mathfrak r}

\def\k{{\boldsymbol k}}

\def\wi{\widetilde}		
	\def\deg{\operatorname{deg}}		
\def\Ind{\operatorname{Ind}}	
\def\Gal{\operatorname{Gal}}	\def\Rad{\operatorname{Rad}}
\def\Ker{\operatorname{Ker}}	
 \def\Im{\operatorname{Im}}

\date{}

\title{\LARGE{Principalisation des groupes de classes logarithmiques}}

\author{Jean-François {\sc Jaulent}}

\begin{document}
\maketitle
\medskip

{\footnotesize
\noindent{\bf Résumé.} 
Nous transposons aux classes logarithmiques attachées à un corps de nombres les résultats sur la principalisation abélienne des groupes de classes de rayons modérées obtenus dans un article antérieur.
}\smallskip

{\footnotesize
\noindent{\bf Abstract.} 
We extend to logarithmic class groups the results on abelian principalization of tame ray class groups of a number field obtained in a previous article
}
\smallskip\bigskip

\section*{Introduction}
\addcontentsline{toc}{section}{Introduction}
\medskip

Le $\ell$-groupe des classes logarithmiques (de degré nul) $\,\wCl_K$ d'un corps de nombres $K$ a été introduit dans \cite{J28} et se présente comme un analogue formel du $\ell$-sous-groupe de Sylow du groupe $Cl_K$ des classes d'idéaux de ce corps pour un premier $\ell$ donné. Son calcul effectif a été récemment implanté par Belabas dans le système {\sc pari}  (cf. \cite{BJ}). \smallskip

Par la Théorie $\ell$-adique du corps de classes (cf. e.g. \cite{J31}), le groupe $\,\wCl_K$ s'interprète comme groupe de Galois $\Gal(K^{lc}/K^c)$ de la pro-$\ell$-extension abélienne localement cyclotomique maximale $K^{lc}$ de $K$ relativement à la $\Zl$-extension cyclotomique $K^c$. En d'autres termes, le groupe $\,\wCl_K$ mesure l'écart pour une $\ell$-extension abélienne entre être localement ou globalement cyclotomique. Cela explique le rôle souvent implicite que jouent ces groupes logarithmiques dans l'étude des pro-$\ell$-extensions cyclotomiques, notamment dans l'interprétation de la conjecture de Greenberg évoquée en appendice du présent travail.\smallskip

Du point de vue local, le passage de la valuation classique à la valuation logarithmique revient à remplacer la $\Zl$-extension non-ramifiée par la $\Zl$-extension cyclotomique, ce qui permet de définir les notions de degré d'inertie et d'indice de ramification au sens logarithmique en analogie avec les mêmes objets traditionnels. Les extension logarithmiquement non-ramifiées sont ainsi les extensions localement cyclotomiques. De plus, les unités du corps local au sens logarithmique sont tout simplement les normes cyclotomiques locales. De ce fait, les unités logarithmiques globales sont exactement les normes cyclotomiques. Leur rang est donné par la conjecture de Gross-Kuz'min (initialement énoncée dans \cite{Kuz}), qui revient à postuler qu'il est égal à celui des unités ordinaires augmenté de 1 (cf. \cite{J28,J31,J55}) ou, de façon équivalente, que le pro-$\ell$-groupe $\,\wCl_K$ est fini. C'est, en particulier le cas lorsque le corps $K$ est abélien. Il en est de même pour certaines familles de corps non abéliens, dits $\ell$-rationnels, pour lesquelles on peut montrer qu'il est trivial (cf. \cite{GJ,JN,Mov}).

En résumé, les $\ell$-groupes de classes logarithmiques se comportent comme les $\ell$-sous-groupes de Sylow des groupes de classes habituels (ce qui permet par exemple de construire des $\ell$-tours localement cyclotomiques analogues aux $\ell$-tours de corps de classes de Hilbert (cf. \cite{JM1,JM2,JS}), avec cependant des différences essentielles: en particulier le théorème 94 de Hilbert, qui joue un rôle clé dans les questions de capitulation, ne s'applique pas dans le cadre logarithmique (cf. \cite{J54}).\medskip

L'objet du présent article est ainsi de reprendre dans le cadre logarithmique les travaux de \cite{J61} qui généralisent aux classes de rayons les résultats antérieurs de Gras \cite{Gr1,Gr6}, Kurihara \cite{Kuh} et Bosca \cite{Bo1,Bo2} sur la principalisation abélienne des groupes de classes d'idéaux. Nous avons fait le choix de suivre aussi fidèlement que possible la démarche de \cite{J61}, pour faciliter la comparaison et mettre en relief similitudes et dissemblances d'avec le cas classique, la principale étant la nécessité de se restreindre aux classes relatives du fait de la formule du produit (ou du degré) qui n'intervient pas dans le cas des idéaux, mais vient ici compliquer la démonstration.

\newpage
\section{Complément sur les classes logarithmiques de rayons}

Classiquement, le groupe des classes d'idéaux d'un corps de nombres $K$ est défini comme conoyau $Cl_K^{\scale{.8}{\mathrm{ord}}}$ du morphisme naturel partant du groupe multiplicatif $K^\times$ à valeurs dans le groupe des idéaux $I_K$ donné par la famille des valuations $\nu=(\nu_\p)_{\p\in Pl_K^\circ}$ attachées aux places finies de $K$.\smallskip

\centerline{$1 \rightarrow E_K \rightarrow  K^\times \overset{\nu}{\longrightarrow}  I_K \rightarrow Cl^{\scale{.8}{\mathrm{ord}}}_K \rightarrow 1$.}\smallskip

\noindent Fixons maintenant un nombre premier $\ell$. Par produit tensoriel avec $\Zl$, le $\ell$-sous-groupe de Sylow  $\Cl^{\scale{.8}{ord}}_K$ de $Cl^{\scale{.8}{ord}}_K$ apparaît alors comme conoyau du morphisme $\nu$ étendu au tensorisé $\R_K=\Zl\otimes_\ZZ K^\times$ et à valeurs dans le $\Zl$-module libre construit sur ces mêmes places: $\Dl_K=\oplus_{\p\in Pl_K^\circ}\Zl\p$.
\smallskip

Le $\ell$-groupe des classes logarithmiques est le groupe analogue $\,\wCl_K$ obtenu en remplaçant les valuations classiques $\nu_\p$ par leurs homologues $\ell$-adiques $\wi\nu_\p$ définis à partir des logarithmes des valeurs absolues $\ell$-adiques et en se restreignant aux diviseurs de degré nul (cf. \cite{J28, J31}):

\centerline{$1 \rightarrow \wE_K \rightarrow  \R_K\overset{\wi\nu}{\longrightarrow}  \wDl_K \rightarrow \wCl_K \rightarrow 1$.}\smallskip

Pour chaque place finie $\p$ de $K$, notons comme plus haut  $\R_{K_\p}=\varprojlim K_\p^\times/K_\p^{\times\ell^n}$ le compactifié $\ell$-adique du groupe multiplicatif $K_\p^\times$ et $\J_K= \prod_\p^{r\!e\!s} \R_{K_\p}$ le $\ell$-adifié du groupe des idèles de $K$. Introduisons enfin la $\Zl$-extension cyclotomique $K^c$ de $K$.\smallskip

Du point de vue global, la surjection canonique du $\ell$-adifié $\J_K$ du groupe des idèles de $K$ dans le groupe procyclique $\Gal(K^c/K)\simeq\Zl$ fournit  un épimorphisme {\em degré}:

\centerline{$\deg:\J_K\rightarrow \Zl$;}

\noindent dont le noyau $\wi\J_K$ est, par construction, le sous-groupe normique de $\J_K$ attaché à $K^c$. Son quotient\smallskip

\centerline{$\,\wDl_K=\wJ_K/\wU_K$}\smallskip

\noindent par le sous-groupe $\,\wU_K=\prod_{\p}\wU_\p$ est le {\em $\ell$-groupe des diviseurs logarithmiques de degré nul}. L'image\smallskip

\centerline{$\,\Pl_K\simeq\R_K/\wE_K$}\smallskip

\noindent de $\,\R_K$ dans $\,\wDl_k$ est le {\em sous-groupe des diviseurs logarithmiques principaux}. Et le quotient\smallskip

\centerline{$\,\wCl_K = \wDl_K/\Pl_K \simeq \wJ_K/\prod_\p\wU_\p\R_K$}\smallskip

\noindent est, par construction, le {\em $\ell$-groupe des classes logarithmiques du corps $K$}.
La {\em conjecture de Gross-Kuz'min} pour le corps $K$ et le premier $\ell$ en postule la finitude (cf. \cite{J28,J31,J55}).\smallskip

 Du point de vue local, le noyau $\,\U_\p$ de $\nu_\p$ dans $\R_\p$ (autrement dit le sous-groupe des unités de $\R_\p$) est le groupe de normes associé à la $\Zl$-extension non ramifiée de $K_\p$; tandis que le noyau $\,\wU_\p$ de $\wi\nu_\p$ (i.e. le sous-groupe des unités logarithmiques) correspond, lui, à sa $\Zl$-extension cyclotomique. Par analogie avec le cas classique, il est commode de dire qu'une $\ell$-extension localement cyclotomique est {\em logarithmiquement non-ramifiée}.\smallskip
 
Par la Théorie $\ell$-adique du corps de classes (cf. \cite{Gr2,J28,J31}), le $\ell$-groupe des classes logarithmiques d'idéaux s'interprète comme groupe de Galois $\Gal(K^{lc}/K^c)$ de la pro-$\ell$-extension abélienne localement cyclotomique maximale $K^{lc}$ de $K$ relativement à la $\Zl$-extension cyclotomique $K^c$. En  particulier $K^{lc}$ est la plus grande pro-$\ell$-extension abélienne de $K$ qui est complètement décomposée au-dessus de $K^c$, i.e. logarithmiquement non-ramifiée sur $K$. Plus généralement:

\begin{Def}
Étant donné un ensemble fini $T_K$ d'idéaux premiers de $K$ étrangers à $\ell$ et $\m_K^\ph$ le produit $\prod_{\q_K^\ph\in T_K^\ph}\q_K^\ph$, le $\ell$-groupe des classes logarithmiques de rayons  modulo $\m$ est le quotient

\centerline{$\wCl^{\,\m_{\si{K}}}_K=\wDl^{\,\m_{\si{K}}}_K/\Pl^{\,\m_{\si{K}}}_K$}\smallskip

\noindent du $\Zl$-module $\wDl^{\,\m_{\si{K}}}_K$ construit sur les diviseurs logarithmiques étrangers à $T_K$ par l'image $\Pl^{\,\m_{\si{K}}}_K$ du sous-module $T_K$-infinitésimal $\R^{\,\m_{\si{K}}}_K=\big\{x\in\R_K\;|\; s_{\q_{\si{K}}}(x)=1  \; \forall\q_K^\ph\in T_K^\ph \big\}$ de  $\R_K=\Zl\otimes_\ZZ K^\times$.
\end{Def}

\begin{Sco}
Le $\ell$-groupe $\,\wCl^{\,\m}_K$ s'interprète comme le groupe de Galois $\Gal(\wi H_K^T/K^c)$ attaché à  l'extension abélienne $T$-logarithmiquement ramifiée (i.e. non-ramifiée au sens logarithmique en dehors de $T$) maximale $\wi H_K^{T_{\si{K}}}$ attachée à $K$ relativement à la $\Zl$-extension cyclotomique $K^c$.
\end{Sco}

\Preuve Par un calcul immédiat, on a, en effet: $\,\wCl^{\,\m_{\si{K}}}_K = \wDl^{\,\m_{\si{K}}}_K/\Pl^{\,\m_{\si{K}}}_K \simeq \wJ_K/\prod_{\q_K^\ph\notin T_K^\ph}\wU_{\q_{\si{K}}^\ph}\R_K^\ph$.

\newpage
\section{Énoncé du résultat et stratégie de preuve}

La situation considérée est la suivante: $\ell$ est un nombre premier fixé; $K/\k$ désigne une extension galoisienne de corps de nombres; et $T=T_\k$ est un ensemble fini de places finies de $\k$ ne divisant pas $\ell$. Pour chaque extension finie $N$  de $\k$, nous notons $T_N$ l'ensemble des premiers de $N$ au-dessus de $T$ et $\m^\ph_N=\prod_{\q_N^\ph\in T_N^\ph}\q^\ph_N$ leur produit.\smallskip

Nous nous proposons de faire capituler le $\ell$-groupe des classes logarithmiques de rayons attachées à $K$ par composition avec une $\ell$-extension {\em abélienne} $F$ de $\k$; plus précisément de prouver que $\,\Cl^{\,\m_{\si{K}}}_K$ a une image triviale dans  $\,\Cl^{\,\m_{\si{L}}}_L$ pour $L=KF$ et une infinité de telles extensions $F$.

Néanmoins, pour des raisons spécifiques aux classes logarithmiques, il est naturel pour cela de restreindre notre ambition aux classes relatives:

\begin{Def}\label{CLR}
Par sous-groupe des classes relatives du $\ell$-groupe des classes logarithmiques de rayons modulo $\m_K^\ph$ nous entendons le noyau $\,\Cl^{\,\m_{\si{K}}}_{K/\k}$  de l'application norme $N_{K/\k}:\Cl^{\,\m_{\si{K}}}_K\to\Cl^{\,\m_\k^\ph}_\k$.
\end{Def}

Rappelons enfin que la conjecture de Gross-Kuz'min pour le premier $\ell$ et le corps $N$ revient à postuler la finitude du $\ell$-groupe des classes logarithmiques $\,\wCl_N$ attaché à ce corps.

Cela étant, le résultat principal de cette note peut s'énoncer comme suit:

\begin{Th}\label{ThP}
Soient $K$ un corps de nombres qui satisfait la conjecture de Gross-Kuz'min pour un premier $\ell$ et $\k$ un sous-corps tel que $K/\k$ soit complètement décomposé en au moins une place à l'infini.
Pour tout ensemble fini $T$ de places de $\k$ ne divisant pas $\ell$, il existe une infinité de $\ell$-extensions abéliennes  $F/\k$ complètement décomposées en toutes les places à l'infini, telles que le sous-groupe relatif $\,\Cl^{\,\m_{\si{K}}}_{K/\k}$ du $\ell$-groupe $\,\Cl^{\,\m_{\si{K}}}_K$ des classes logarithmiques de rayons modulo $\m_K^\ph=\prod_{\q_K^\ph \in T_K^\ph} \q_K^\ph$  capitule dans le compositum $L=FK$.
\end{Th}

\begin{Cor}
Sous les hypothèses du Théorème, dès lors que le $\ell$-groupe $\,\wCl^{\,\m_\k^\ph}_\k$ des classes logarithmiques de rayons du corps de base $\k$ est trivial, le $\ell$-groupe $\,\wCl^{\,\m_{\si{K}}}_K$ entier capitule dans $FK$.
\end{Cor}

Venons-en maintenant à la stratégie de la preuve. Elle est essentiellement analogue à celle utilisée par Bosca pour principaliser les classes d'idéaux (cf. \cite{Bo1,Bo2}) et récemment étendue aux classes de rayons classiques (cf. \cite{J61}) avec toutefois quelques complications, la première étant que le pro-$\ell$-groupe des classes logarithmiques sans restriction de degré, disons $\,\Cl^{\si{\,log}}_K \simeq \Gal(K^{lc}/K)$, est infini: dans le cas des groupes classes d'idéaux ou de rayons, qui sont finis, il est toujours possible de représenter une classe donnée par un idéal premier (auquel on impose des conditions supplémentaires ad hoc). Mais c'est impossible dans le cas logarithmique, les diviseurs premiers n'étant jamais de degré nul. Il faut donc biaiser.\smallskip

Étant donnée une classe (de degré nul) $[\,\wi\d_K^\ph]$ dans $\,\wCl^{\,\m_{\si{K}}}_K$, un entier $n$ (assez grand) ayant été choisi, nous pouvons cependant écrire

\centerline{$[\,\wi\d_K^\ph]=[\,\wi\p_K^\ph]+\ell^n [\,\wi\b_K^\ph]$,}\smallskip

\noindent avec $\p_K^\ph$ premier pour un diviseur logarithmique convenable $\wi\b_K^\ph$. Notant alors $\p$ l'unique diviseur premier de $\k$ au-dessous de $\p_K^\ph$, puis imposant à la $\ell$-extension cyclique $F_\d/\k$ d'être logarithmiquement non-ramifiée en dehors de $\p$ et d'avoir pour indice de ramification logarithmique $\wi e_{\p} (F_\d/\k)=\ell^n$, nous obtenons dans le pro-$\ell$-groupe des diviseurs logarithmiques du compositum $L_\d=F_\d K$ l'identité entre diviseurs: $\wi\p_K^\ph=\ell^n\; \wi\a_{L_\d}^\ph$, pour un certain diviseur ambige $\wi\a_{L_\d}^\ph$, ainsi dénommé car invariant par $\Gal(F_\d/\k)$. En fin de compte, il vient:\smallskip

\centerline{$[\,\wi\d_K^\ph]=\ell^n[\, \wi\a_{L_\d}^\ph]+\ell^n [\,\wi\b_K^\ph] = \ell^n [\, \wi\a_{L_\d}^\ph +\,\wi\b_K^\ph ]$,}\smallskip

\noindent dans $\,\wCl^{\,\m_{\si{L_\d}}^\ph}_{L_\d}$, pour un certain diviseur ambige $\wi\a_{L_\d}^\ph +\,\wi\b_K^\ph$; de sorte que $\wi\d_K^\ph$ se principalise dans $L_\d$, dès lors que $\ell^n$ annule le sous-groupe ambige qui en contient la classe.\smallskip

Finalement, prenant le compositum $F$ des $F_\d$ pour un système de représentants de générateurs de $\,\wCl^{\,\m_{\si{K}}}_K$ et posant $L=FK$, on obtient bien un corps principalisant pour $\,\wCl^{\,\m_{\si{K}}}_K$.\smallskip

La première étape consiste donc à préciser le nombre de classes d'ambiges dans une $\ell$-extension cyclique, pour pouvoir le majorer indépendamment de $n$ sous certaines conditions.

\newpage
\section{La formule des classes logarithmiques d'ambiges}

Nous reprenons ci-dessous, en les modifiant légèrement pour les adapter aux classes logarithmiques de rayons, les calculs de classes invariantes effectués dans \cite{J28}.

Supposons donc fixés un nombre premier $\ell$ et une $\ell$-extension cyclique de corps de nombres $L/K$; donnons-nous un ensemble fini $T_K^\ph$ d'idéaux premiers de $K$ ne divisant pas $\ell$ notons $T_L^\ph$ l'ensemble des premiers de $L$ au-dessus de $T_K^{\si{(\ell)}}$; posons enfin $\m^\ph_K=\prod_{\q^\ph_K\in T_K^\ph}\q_K$ et $\m^\ph_L=\prod_{\q^\ph_L \in T_L^\ph}\q^\ph_L$.
Le résultat logarithmique s'énonce alors, en analogie avec le résultat classique (cf. \cite{J61}, Prop. 3):

\begin{Prop}\label{CLA}
Dans une $\ell$-extension cyclique $L/K$ de corps de nombres, le nombre de classes logarithmiques de rayons dans $\wCl_L^{\,\m_{\si{L}}}$ qui sont représentées par des diviseurs logarithmiques  ambiges (i.e. invariants par $C=\Gal(L/K)$) est donné sous la conjecture de Gross-Kuz'min par la formule:
$$
(\wDl_L^{\,\m_{\si{L}}}{}^C:\Pl_L^{\,\m_{\si{L}}}{}^C) \;=\; |\,\wCl_K^{\,\m_{\si{K}}}|\;\frac{\prod_{\p_{\si{\!\infty}}} d_{\p_{\si{\!\infty}}}\!(L/K)\,\prod_{\p_{\si{\circ}}\nmid\m} \tilde e_{\p_{\si{\circ}}}(L/K)}{\big(\deg_L\wDl_L^{\,\m_{\si{L}}}{}^C:\deg_L\wDl^{\,\m_{\si{K}}}_K\big)\;\big( \wE^{\,\m_{\si{K}}}_K:N_{L/K}(\wE^{\,\m_{\si{L}}}_L) \big)}
$$ 
\noindent Dans celle-ci $\,d_{\p_{\si{\!\infty}}}\!$ désigne le degré local et $\tilde e_{\p_{\si{\circ}}}$  l'indice de ramification logarithmique; $\p_{\si{\!\infty}}$ parcourt les places à l'infini de $K$ et $\p_{\si{\circ}}$ les places étrangères à $\m^\ph_K$ logarithmiquement ramifiées dans $L/K$.
\end{Prop}

\Preuve Écrivons: $(\wDl_L^{\,\m_{\si{L}}}{}^C:\Pl_L^{\,\m_{\si{L}}}{}^C) \;=\; (\wDl_L^{\,\m_{\si{L}}}{}^C:\wDl^{\,\m_{\si{K}}}_K)\,(\wDl^{\,\m_{\si{K}}}_K:\Pl^{\,\m_{\si{K}}}_K)\,/\,(\Pl_L^{\,\m_{\si{L}}}{}^C:\Pl^{\,\m_{\si{K}}}_K)$.\smallskip

\begin{itemize}
\item  L'indice $ (\wDl_L^{\,\m_{\si{L}}}{}^C:\wDl^{\,\m_{\si{K}}}_K)$ se calcule comme suit: dans la suite exacte courte canonique\smallskip

\centerline{$1 \rightarrow \wDl_L^{\,\m_{\si{L}}}{}^C / \wDl^{\,\m_{\si{K}}}_K \rightarrow \Dl_L^{\,\m_{\si{L}}}{}^C / \Dl^{\,\m_{\si{K}}}_K \rightarrow \Dl_L^{\,\m_{\si{L}}}{}^C/\wDl_L^{\,\m_{\si{L}}}{}^C\Dl^{\,\m_{\si{K}}}_K \rightarrow 1$,}\smallskip

\noindent le terme de droite s'identifie via l'application {\em degré} au quotient $\deg_L\wDl_L^{\,\m_{\si{L}}}{}^C/\deg_L\wDl^{\,\m_{\si{K}}}_K$ .
Et le groupe $\Dl_L^{\,\m_{\si{L}}}{}^C$ des diviseurs logarithmiques ambiges étrangers à $\m^\ph_L$ est engendré par les sommes $\sum_{\p^\ph_L|\p^\ph_K}\p^\ph_L=\tilde e^{\si{-1}}_{\p^\ph_{\si{K}}}(L/K)\,\p^\ph_K$, lorsque $\p^\ph_K$ décrit l'ensemble des premiers de $K$ qui ne divisent pas $\m_K^\ph$, d'où:
 
\centerline{$(\Dl_L^{\,\m_{\si{L}}}{}^C:\Dl^{\,\m_{\si{K}}}_K)=\prod_{\p_{\si{\circ}}\nmid\m_K^\ph}\tilde e_\p(L/K)$.}\smallskip

\item Le quotient $\wDl^{\,\m_{\si{K}}}_K/\Pl^{\,\m_{\si{K}}}_K$ est tout simplement le $\ell$-groupe $\,\wCl^{\,\m_{\si{K}}}_K$ des classes logarithmiques de $K$, lequel est fini sous la conjecture de Gross-Kuz'min dans $K$.\smallskip

\item Enfin, tout comme dans le cas classique (cf. \cite{J61}), le lemme du serpent, appliqué ici à la suite exacte de cohomologie associée à la suite courte qui définit le sous-groupe principal\smallskip

\centerline{$1 \longrightarrow \wE^{\,\m_{\si{L}}}_L \longrightarrow \R^{\m_{\si{L}}}_L \longrightarrow \Pl^{\,\m_{\si{L}}}_L \longrightarrow 1$,}\smallskip

montre que le quotient $\Pl_L^{\,\m_{\si{L}}}{}^C/\Pl^{\,\m_{\si{K}}}_K$ s'identifie au premier groupe de cohomologie des unités logarithmiques $T_L^{\si{(\ell)}}$-infinitésimales $H^1(C,\wE^{\,\m_{\si{L}}}_L)$:\smallskip

\centerline{$\Pl_L^{\,\m_{\si{L}}}{}^C/\Pl^{\,\m_{\si{K}}}_K \simeq  H^1(C,\wE^{\,\m_{\si{L}}}_L)$,}\smallskip

\noindent dont l'ordre est le produit de celui du groupe $\wE^{\,\m_{\si{K}}}_K/N_{L/K}(\wE^{\,\m_{\si{L}}}_L)) \simeq H^2(C,\wE^{\,\m_{\si{L}}}_L)$ par le quotient de Herbrand $q(C,\wE^{\,\m_{\si{L}}}_L)$, lequel ne dépend que du caractère des unités logarithmiques. Sous la conjecture de Gross-Kuz'min, ce dernier est donné par la formule (cf. \cite{J28}, Th. 3.6):\smallskip

\centerline{$\chi_{\wE^{\m_{\si{L}}^\ph}_L}=\chi^\ph_{\wE^\ph_L}=\sum_{\p_{\si{\!\infty}}} \Ind^C_{D_{\p_{\si{\!\infty}}}} 1_{D_{\p}}$,}

\noindent comme somme d'induits attachés aux sous-groupes de décomposition des places à l'infini de $K$. Il vient donc ici:

\centerline{$q(C, \wE_L^{\,\m_{\si{L}}})=q(C,\wE_L)=\prod_{\p_{\si{\!\infty}}}d_{\p_{\si{\!\infty}}}\!(L/K)$.}\smallskip

\noindent D'où la formule annoncée.
\end{itemize}\medskip

\Remarque Dans l'expression obtenue, le facteur $\big(\deg_L\wDl_L^{\,\m_{\si{L}}}{}^C:\deg_L\wDl^{\,\m_{\si{K}}}_K\big)$ vient remplacer le facteur $[L:K]$ au dénominateur de la formule analogue pour les classes de rayons ordinaires.
Lorsque l'extension $L/K$ est totalement ramifiée au sens logarithmique en un diviseur {\em primitif} $\,\p_K$ (i.e. tel que $\deg_K\p^\ph_K$ engendre $\Gal(K^c/K)$) qui ne divise pas $\m_K$, il y a même égalité:\smallskip

\centerline{$\big(\deg_L\wDl_L^{\,\m_{\si{L}}}{}^C:\deg_L\wDl^{\,\m_{\si{K}}}_K\big))=[L^c:K^c]=[L:K]$.}\smallskip

\noindent En général, cependant, on a simplement les relations de divisibilité.

\newpage
\section{Classes relatives d'ambiges et modules équivalents}

Revenons maintenant au problème de la principalisation. En analogie avec le cas des classes de rayons classiques, nous cherchons  $L$ comme compositum $FK$ pour une $\ell$-extension cyclique $F/\k$ logarithmiquement ramifiée en une unique place $\p$ de $\k$, complètement décomposée dans $K/\k$, ayant par ailleurs un indice de ramification logarithmique $\wi e_{\p}(F/\k)=\ell^n$ suffisamment grand.

Malheureusement la condition de décomposition requise, en tuant le facteur degré au dénominateur, ruine toute possibilité de contrôler le quotient indépendamment de $n$ en toute généralité.

Pour pallier cette difficulté, nous allons nous restreindre au sous-groupe des classes relatives, i.e. au noyau de la norme arithmétique $N_{K/\k}$ attachée à l'extension $K/\k$ (cf. Définition \ref{CLR}).\medskip

Notons $d$ l'ordre de $G=\Gal(K/\k)$; puis $\nu=\sum_{\sigma\in G}\sigma$ l'élément de $G$ qui correspond à la norme $N_{K/\k}$; et $\bar\nu=\sum_{\sigma\in G}(1-\sigma) =d-\nu$. Le noyau de la norme {\em arithmétique}  $N_{K/\k}$ est évidemment compris entre le noyau de la norme {\em algébrique} $\nu$  et l'image de l'opérateur complémentaire $\bar\nu$:

\centerline{$\Im\bar\nu \subset \Ker N_{K/\k} \subset \Ker \nu$.}
\begin{itemize}
\item Lorsque l'ordre $d$ de $G$ est étranger à $\ell$, il est inversible dans l'anneau $\Zl$, de sorte que  les éléments $e_\nu=\frac{\nu}{d}$ et $e_{\bar\nu}=\frac{\bar\nu}{d}$ sont deux idempotents centraux complémentaires de l'algèbre $\Zl[G]$. Et il vient donc: $\Im\bar\nu=\Im e_{\bar\nu} =\Ker e_\nu =\Ker \nu$. Dans ce cas, il n'est pas nécessaire de distinguer entre noyau arithmétique et noyau algébrique de la norme. Mieux encore, toute suite exacte de $\Zl[G]$-modules donne aussitôt deux suites exactes, l'une restreinte aux noyaux de la norme $\nu$, l'autre aux images, par action des deux idempotents précédents.
\item Dans le cas contraire, il faut distinguer. Cependant, eu égard au problème qui nous préoccupe, il est possible de contourner cette difficulté supplémentaire de la façon suivante:
\end{itemize}

\begin{Conv}
Convenons de dire que deux $\Zl[G]$ modules finis $M$ et $N$ dépendant du paramètre $n$ sont équivalents lorsque la $\ell$-valuation du quotient de leurs ordres est bornée indépendamment de $n$; ce que nous écrivons:

\centerline{$M \sim N \iff |M| \approx |N] \iff v_\ell\big(\frac{|M|}{|N|}\big)$ borné (indépendamment de $n$).}
\end{Conv}


\begin{Lem}
Pour tout $\Zl[G]$-module fini $M$ de $\Zl$-rang borné (indépendamment de $n$), le noyau ${}_\nu M$ de $\nu$ et l'image $M^{\bar\nu}$ de $\bar\nu$ sont équivalents; de même le noyau ${}_{\bar\nu} M$ de $\bar\nu$ et l'image $M^{\nu}$ de $\nu$.
\end{Lem}

\Preuve L'identité $\nu+\bar\nu=d$ montre que les deux quotients ${}_\nu M/M^{\bar\nu}$ et ${}_{\bar\nu} M/M^{\nu}$ sont tués par $d$.\medskip

En conséquence, quitte à travailler à un borné près, il est toujours possible, dès lors qu'on  ne considère que des modules de rang borné (indépendamment de $n$), de raisonner comme si la condition de semi-simplicité $\ell\nmid d$ était toujours remplie. En particulier, pour toute suite exacte\smallskip

\centerline{$1 \to N \to M \to P \to 1$}\smallskip

\noindent de  $\Zl[G]$-modules finis de $\Zl$-rang borné  (indépendamment de $n$), on a les équivalences:\smallskip

\centerline{$M^{\bar\nu} \sim {}_\nu M \sim {}_\nu N \oplus {}_\nu P \sim N^{\bar\nu} \oplus P^{\bar\nu} \qquad \& \qquad 
M^{\nu} \sim {}_{\bar\nu} M \sim {}_{\bar\nu} N \oplus {}_{\bar\nu} P \sim N^{\nu} \oplus P^{\nu}$.}\smallskip

\noindent entre pseudo-composantes {\em relatives} (i.e. tuées par $\nu$) à gauche et {\em induites}  (tuées par $\bar\nu$) à droite.\medskip

Reprenons maintenant les calculs de classes d'ambiges effectués dans la section précédente. Dans la situation galoisienne considérée ici, les divers $\ell$-groupes qui interviennent dans  la preuve de la Proposition \ref{CLA} sont des $\Gal(L/\k)$-modules fixés par $C=\Gal(L/K)$, i.e. des $\Zl[G]$-modules finis.
En nous restreignant aux pseudo-composantes relatives (i.e. tuées par $\nu=\sum_{\sigma\in G}\sigma$) et en négligeant tout les modules d'ordre borné indépendamment de $n$, nous obtenons immédiatement:

\begin{Prop}\label{CAR}
Sous les hypothèses énoncées en début de section, l'ordre de la pseudo-composante relative ${}_\nu\big(\wDl_L^{\,\m_{\si{L}}}{}^C/\Pl_L^{\,\m_{\si{L}}}{}^C\big)$ du sous-groupe des classes logarithmiques de rayons dans $\wCl_L^{\,\m_{\si{L}}}$ qui sont représentées par des diviseurs logarithmiques  ambiges (i.e. invariants par $C=\Gal(L/K)$) est donné sous la conjecture de Gross-Kuz'min (pour le premier $\ell$ et le corps $L$)  par la formule:\smallskip

\centerline{$\big|\,{}_\nu\big(\wDl_L^{\,\m_{\si{L}}}{}^C/\Pl_L^{\,\m_{\si{L}}}{}^C\big)\,\big| \approx \ell^{n(d-1)} / \big ({}_\nu\wE^{\,\m_{\si{K}}}_K /({}_\nu\wE^{\,\m_{\si{K}}}_K \cap N_{L/K}(\wE^{\,\m_{\si{L}}}_L))\big)$.}
\end{Prop}

\Preuve De l'isomorphisme $\Dl_L^{\,\m_{\si{L}}}{}^C/\Dl^{\,\m_{\si{K}}}_K\simeq(\ZZ/\ell^n\ZZ)[G]$, on tire: $|\,{}_\nu (\Dl_L^{\,\m_{\si{L}}}{}^C/\Dl^{\,\m_{\si{K}}}_K)\,|=\ell^{n(d-1)}$.

\newpage
\section{Minoration de l'indice normique des unités}

Sous la conjecture de Gross-Kuz'min (pour le premier $\ell$ et le corps $L$) le théorème de représentation des unités logarithmiques rappelé plus haut (cf. \cite{J28}, Th. 3.6) appliqué à l'extension galoisienne $L/K$ nous assure que $\,\wE_K$ et donc son sous-module d'indice fini $\,\wE_K^{\,\m_{\si{K}}}$ contiennent un sous-module monogène $\,\wi\E^{\,\varepsilon}_K=\,\varepsilon^{\Zl[G]}$ de caractère $\chi_G^{\textrm{aug}}$.\smallskip

Rappelons que $\,\wE_K$ est contenu dans le $\ell$-adifié $\,\E'_K=\Zl\otimes_\ZZ E'_K$ du groupe des $\ell$-unités de $K$ (mais qu'il n'est pas en général, comme l'est $\,\E_K=\Zl\otimes_\ZZ E_K$, le $\ell$-adifié d'un sous-groupe de $E'_K$).
Prenons $n$ assez grand; choisissons dans le groupe $E'_K$ un système de représentants $R_K^{\,\varepsilon}$ du quotient $\,\wE_K^{\,\varepsilon} /\big( \wE_K^{\,\varepsilon}\cap\E'_K{\!}^{\ell^n} \big)$; notons $\mu_{\ell^n}$ le $\ell$-groupe des racines $\ell^n$-ièmes de l'unité;
 et considérons l'extension $ K_n^{\varepsilon}=K\big[\mu^\ph_{\ell^n},\sqrt[\ell^n]{R^{\,\varepsilon}_K}\,\big]$ engendrée sur $K_n=K[\mu^\ph_{\ell^n}]$ par les racines $\ell^n$-ièmes des éléments de $R^{\,\varepsilon}_K$.  
 Notons enfin $\,\wE^{\,\circ}_{K_n}$ la racine de $\,\wE^{\,\varepsilon}_K$ dans $\,\wE^\ph_{K_n}$, puis $\ell^{\tilde\delta_\m}$ l'indice $(\wE^{\,\circ}_{K_n}\!:\wE^{\,\varepsilon}_K\mu^\ph_{K_n})$, qui est ultimement indépendant de $n$. Par construction, le degré de l'extension kummérienne  $K_n^\varepsilon/K_n$ est alors donné par la  formule:\smallskip

\centerline{$[K_n^\varepsilon : K_n]\,=\,|\Rad(K_n^\varepsilon/K_n)|=\big(\wE^{\,\varepsilon}_K : \wE^{\,\varepsilon}_K{}^{\ell^n}\big) /\big(\wE^{\,\circ}_{K_n}:\wE^{\,\varepsilon}_K\mu^\ph_{K_n}\big) \,=\, \ell^{n(d-1)}/\ell^{\tilde\delta_\m}$.}\medskip

Tout comme dans  \cite{J61} \S4, introduisons la sous-extension élémentaire $K_n^e/K_n$ de $K_n^\varepsilon/ K_n$; notons $(\eta_1K_n^{\times\ell},\cdots,\eta_tK_n^{\times\ell})$
 une $\Fl$-base de $\Rad(K_n^e/K_n)\subset K_n^\times/K_n^{\times\ell}$ représentée par des conjugués $\eta_i=\varepsilon^{\sigma_i}$ de $\varepsilon$ dont $\eta_1=\varepsilon$; 
 et définissons $\bar\tau_n \in \Gal(K_n^e/K_n)$ par: $\sqrt[\ell]{\eta_1}^{\,(\bar\tau_n-1)}=\zeta_\ell$
  \quad \& \quad $\sqrt[\ell]{\eta_i}^{\,(\bar\tau_n-1)}=1$, pour $i=2,\cdots,t$.
 Faisons choix enfin d'un relèvement $\tau_n$ de $\bar\tau_n$ dans $G_n$. Il est clair que les conjugués de $\bar\tau_n$ engendrent $G_n/G_n^\ell= \Gal(K_n^e/K_n)$ donc que les conjugués de $\tau_n$ engendrent $G_n$.\smallskip

Supposons choisi enfin un entier $m=m_K$ indépendant de $n$ et soit alors une place $\p^\ph_{\si{K}} \nmid 2\ell\m^\ph_{\si{K}}$ de $K$ au-dessus d'un premier $p$, complètement décomposé dans $K_n/\QQ$, telle que l'application de Frobenius associée à l'extension abélienne $K_n^{\,\varepsilon}/K_n$ envoie l'une des places $\p^\ph_{\si{K_n}}$ de $K_n$ au-dessus de $\p_{\si{K}}^\ph$ sur la puissance $\ell^m$-ième de l'automorphisme $\tau_n$. Notant $s_{\m_{\si K}}$ l'épimorphisme de semi-localisation de $\R_K$ sur $\R_{K_{\m_{\si K}}}=\prod_{\q_{\si{K}}\mid\m_{\si{K}}}\R_{K_{\q_{\si{K}}}}$, nous concluons à l'identité:\smallskip

\centerline{$\big\{\eta\in\wE_K^{\,\varepsilon} \,|\, s_{\m_{\si K}}\!(\eta)\in  \underset{\sigma\in G}\prod \U_{K_{\p_{\si K}^{\si \sigma}}}^{\ell^n}\big\}
= \big\{\eta\in\wE_K^{\,\varepsilon} \,|\, s_{\m_{\si K}}\!(\eta)\in \underset{\sigma\in G}\prod  \R_{K_{\p_{\si K}^{\si \sigma}}}^{\ell^n}\big\}
= \wE_K^{\,\varepsilon}\cap \R_{K_n}^{\times\ell^{n-m}}$.}\medskip

Si donc la $\ell$-extension cyclique $F/\k$ est logarithmiquement non-ramifiée en dehors de la place $\p=\p_\k$ au-dessous de $\p_K$ avec un indice de ramification logarithmique $\tilde e_\p(F/\k)=\ell^n$, les unités logarithmiques de $K$ qui sont normes dans l'extension cyclique $L/K$ sont exactement celles qui sont des puissances $\ell^n$-ièmes locales aux places au-dessus de $\p$; et il vient:\smallskip

\centerline{$ \big(\wE_K^{\,\varepsilon} : \wE_K^{\,\varepsilon} \cap N_{L/K}(\R_L)\big) =
 \big(\wE_K^{\,\varepsilon} : \wE_K^{\,\varepsilon} \cap {}^{\si{-1}}s_{\m_{\si K}}\big( \underset{\sigma\in G}\prod \U_{K_{\p_{\si K}^{\si \sigma}}}^{\ell^n}\big)\big)
= \big(\,\wE_K^{\,\varepsilon}:\wE_K^{\,\varepsilon}\cap \R_{K_n}^{\times\ell^{n-m}}\big) =  \ell^{(n-m)(d-1)}/\ell^{\tilde\delta_\m}$}

\noindent d'où la minoration (au sens de la divisibilité):\smallskip

\centerline{$\big|\,{}_\nu \wE^\m_K/\big({}_\nu \wE^\m_K\cap N_{L/K}(\wE^\m_L)\big)\,\big| \succ 
\big(\wE^{\,\varepsilon}_K:\wE^{\,\varepsilon}_K\cap N_{L/K}(\R_L)\big) \succ  
\ell^{(n-hm)(d-1)-\tilde\delta_\m}$.}\smallskip
 
Il résulte alors de la Proposition \ref{CAR} que l'ordre du $\ell$-groupe ${}_\nu\big(\wDl_L^{\,\m_{\si{L}}}{}^C/\Pl_L^{\,\m_{\si{L}}}{}^C\big)$ est borné indépendamment de $n$.
 
 En résumé, sous la conjecture de Gross-Kuz'min dans $L$, il vient:

 \begin{Prop}\label{CdA}
 Soient $\ell$ un nombre premier, $K/\k$ une extension galoisiennne de corps de nombres dans laquelle une au moins des places à l'infini est complètement décomposée, $T=T_\k$ un ensemble fini de premiers de $\k$ qui ne divisent pas $\ell$ et $T_K$ l'ensemble des premiers de $K$ au-dessus de $T$; puis $\m_\k=\prod_{\q_\k\in T_\k}\q_\k$ et  $\m_K=\prod_{\q_K\in T_K}\q_K$.\smallskip
 
 Soit $n$ un entier assez grand. Si $\p=\p_k \nmid 2\ell\m_\k$ est un idéal premier de $\k$ au-dessus d'un premier $p$ de $\QQ$ complètement décomposé dans l'extension $K[\zeta_{\ell^n}^\ph]/\QQ$ et $F/\k$ une $\ell$-extension cyclique $\p$-logarithmiquement ramifiée avec pour indice de ramification logarithmique $\wi e_\p(F/\k)=\ell^n$ et tel que l'application de Frobenius envoie l'une des places $\p_K$ au-dessus de $\p$ sur une puissance donnée de l'automorphisme $\tau_n$, l'ordre du sous-groupe relatif de $\,\wCl^{\,\m_{\si{L}}}_L$ attaché au compositum $L=\k F$ qui est engendré par les classes logarithmiques des diviseurs invariants par $C=\Gal(L/K)$ est majoré indépendamment de $n$:\smallskip
 
  \centerline{$(\wDl^{\,\m_{\si{L}}}_L{}^C:\Pl^{\,\m_{\si{L}}}_L{}^C) \prec \ell^{\,\tilde c_{K,\m}}$.}
\end{Prop}

\newpage

\section{Construction de l'extension principalisante}

Résumons: étant donnée une extension galoisienne $K/\k$ de corps de nombres complètement décomposée en au moins une place à l'infini, un nombre premier $\ell$ et un diviseur $\m_K$ de $K$ sans facteur carré, étranger à $\ell$ et stable par $G=\Gal(K/\k)$, ayant fait choix d'un sous-module $E_K^{\,\varepsilon}$ de $E_K^{\m_{\si{K}}}$ de caractère $\chi_G^{\textrm{aug}}$ et d'un entier $m=m_K$ indépendant de $n$, nous avons défini une constante $\wi c_{K,\m_{\si{K}}}$; puis, pour une classe $[\d_K]$ d'ordre $\ell$-primaire dans $\wCl_K^{\,\m_{\si{K}}}$ et $n \ge \wi c_{K,\m_{\si{K}}}$, nous cherchons:\smallskip
 
\begin{itemize}
\item une place $\p$ de $\k$ et une place $\p^\ph_K$ de $K$ au-dessus qui satisfasse les trois conditions suivantes:\smallskip
\begin{itemize}
\item[($i$)] $\p_K^\ph$ est au-dessus d'un premier $p\ne\ell$ de $\QQ$ complètement décomposé dans $K_n=K[\zeta_{\ell^n}]$;\smallskip

\item[($ii$)] $\p_K^\ph$ a même image que $\d_K$ dans le quotient d'exposant $\ell^n$ de $\,\wCl_K^{\,\m_{\si{K}}}$, i.e. a même image dans le groupe de Galois $\Gal(\wi H_n/K)$ de la sous-extension d'exposant $\ell^n$ de $\wi H_K^{\,\m_{\si{K}}}/K$;\smallskip

\item[($iii$)]  l'une des places au-dessus $\p_n$ dans $K_n$ est d'image $\tau_n^{\ell^h}$ dans $\Gal\big(K_n\big[\!\sqrt[\ell^{\si{n}}]{E^{\,\varepsilon}_K}\,\big]/K_n\big)$;
\end{itemize}
\smallskip
\item et une $\ell$-extension $F/\k$ cyclique $\infty$-décomposée, logarithmiquement non-ramifiée en dehors de $\p$ et $\p$-logarithmiquement ramifiée avec pour indice de ramification $\wi e_\p(F/\k)=\ell^n$.\medskip
 \end{itemize}
 
 Examinons d'abord cette dernière condition. Par la théorie $\ell$-adique du corps de classes (cf. \cite{J31}), le groupe de Galois $\Gal(\wi H_k^\p/\wi H^\ph_k)$ de la $\ell$-extension abélienne $\p$-logarithmiquement ramifiée $\infty$-décomposée maximale $\wi H_k^\p$ du corps $k$ relativement à sa sous-extension logarithmiquement non-ramifiée $\infty$-décomposée maximale $\wi H^\ph_k$ est donné par l'isomorphisme:\smallskip
 
 \centerline{$\Gal(\wi H_k^\p/\wi H^\ph_k)\simeq \big(\R_\k\prod_{\q\nmid_\infty}\!\wU_{\k_\q}\prod_{\q\mid_\infty}\!\R_{\k_\q}\big) /
 \big(\R_\k\prod_{\q\ne\p}\wU_{\k_\q}\prod_{\q\mid_\infty}\!\R_{\k_\q}\big) \simeq \mu_{\k_\p}/s_\p(\wE_\k) $,}\smallskip
 
 \noindent où $\mu_{\k_\p}=\,\wU_{\k_\p}$ est le $\ell$-groupe des racines de l'unité dans $\k_\p$ et $s_\p(\wE_\k)$ l'image locale du groupe des unités globales.
Or, le quotient obtenu est cyclique d'ordre $\ell^m$ pour un $m\ge n$ si et seulement si le complété $\k_\p$ contient les racines $\ell^n$-ièmes de l'unité et si les éléments de $\wE_\k$ sont des puissances $\ell^n$-ièmes locales dans $\k_\p$; ce qui a lieu dès que la place $\p$ est complètement décomposée dans l'extension $\k \big[ \zeta_{\ell^n},\wE_\k^{\,\ell^{\scale{0.7}{-n}}}\big]/\k$. Lorsque c'est le cas, l'extension $\wi H_k^\p/\k$ possède donc une sous-extension cyclique $\infty$-décomposée  et $\p$-logarithmiquement ramifiée avec pour indice de ramification $\wi e_\p(F/\k)=\ell^n$.
En fin de compte, l'existence de $F$ est donc assurée dès lors que l'on remplace la condition (i) par:\smallskip
\begin{itemize}
\item[($i'$)]  $\p^\ph_K$ est au-dessus d'un premier $p\ne\ell$ complètement décomposé dans $K_n\big[\wE_\k^{\,\ell^{\scale{0.7}{-n}}}\big]/\QQ$.
\end{itemize}\smallskip

\noindent Et tout le problème est alors de s'assurer de la compatibilité des trois conditions ($i'$), ($ii$) et ($iii$).\smallskip

D'un côté, $\,\wE_\k$ et $\wE_K^{\,\varepsilon}$ étant en somme directe, les extensions $K_n\big[\wE_\k ^{\,\ell^{\scale{0.7}{-n}}}\big]/K_n$ et $K_n\big[ \wE_K^{\,\varepsilon}{}^{\ell^{\scale{0.7}{-n}}}\big]/K_n$ sont bien linéairement disjointes. D'un autre côté, ce n'est pas nécessairement le cas des deux extensions $\wH_K^\m/K_n$ et  $K_n\big[(\wE_\k \wE_K^{\,\varepsilon})^{\ell^{\scale{0.7}{-n}}}\big]/K_n$. Ce défaut demeure cependant borné.\smallskip

Pour voir cela, notons $\ell^{m_{\si{K}}}$ l'ordre de $\mu^\ph_K$; puis $K_{\scale{0.8}{\infty}}=K[\zeta_{\ell^{^{\scale{0.6}{\infty}}}}]$ l'extension cyclotomique  engendrée par toutes les racines d'ordre $\ell$-primaire de l'unité. Observant que $\wi H_K^\m/K$ est abélienne et que $K_n\big[(\wE_\k \wE_K^{\,\varepsilon})^{\ell^{\scale{0.7}{-n}}}\big]/K$ est logarithmiquement non-ramifiée en dehors de $\ell$, nous avons:\smallskip

\centerline{$\wi H_K^\m \cap K_n\big[(\wE^\ph_\k \wE_K^{\,\varepsilon})^{\ell^{\scale{0.7}{-n}}}\big] 
= \wi H_K \cap K_n\big[(\wE^\ph_\k \wE_K^{\,\varepsilon})^{\ell^{\scale{0.7}{-m_{\si{K}}}}}\big]$ pour $n\ge m_{\si{K}}$.}\smallskip

\noindent Et $\Gal(\big[\wi H_K \cap K_n\big[(\wE^\ph_\k \wE_K^{\,\varepsilon})^{\ell^{\scale{0.7}{-m_{\si{K}}}}}\big]:K_n\big]/K_n)$ est un groupe abélien d'exposant divisant $\ell^{m_{\si{K}}}$.
Si donc nous définissons l'entier $m$ introduit dans la section précédente en prenant $m=m_K^\ph$, la compatibilité de $(i')$ avec $(iii)$ est automatique. Celle de $(i')$ avec $(ii)$ est alors immédiate sous la condition supplémentaire:\smallskip
\begin{itemize}
\item[($iv$)]  $[\d_K]$ est une puissance $\ell^{m_{\si{K}}}$-ième dans $\,\wCl_K^{\,\m_{\si{K}}}$.
\end{itemize}\smallskip

Lorsque cette dernière est remplie, le théorème de densité de \v Cebotarev (cf. \cite{Chb} ou  e.g. \cite{Lan}) appliqué dans la clôture galoisienne de l'extension  $\wH_n \big[\zeta_{\ell^n}, (\wE_\k \wE_K^{\,\varepsilon})^{\ell^{\scale{0.7}{-m_{\si{K}}}}}\big]/\QQ$ nous assure l'existence d'une infinité de premiers $p$ possédant une place $\p_K^\ph$ au-dessus qui satisfait les conditions requises $(i')$, $(ii)$ et $(iii)$; et donc l'existence d'une infinité de $\ell$-extensions $F/\k$ convenables. Ainsi:

\begin{Prop}\label{Construction}
La construction de $F$ est possible dès lors que la classe $[\d_K]$ est une puissance $\ell^{m_{\si{K}}}$-ième dans $\,\wCl_K^{\,\m_{\si{K}}}$, où $\ell^{m_{\si{K}}}\!=|\,\mu_K^\ph|$ est l'ordre du $\ell$-groupe des racines de l'unité dans $K$.
\end{Prop}

\newpage
\section{Preuve provisoire du résultat principal}

Plaçons nous d'abord dans le cas galoisien. Sous la conjecture de Gross-Kuz'min, nous avons:

\begin{Th}\label{ThGalois}
Étant donnés une extension galoisienne $K/\k$ de corps de nombres complètement décomposée en au moins une place à l'infini, un nombre premier $\ell$ et un diviseur $\m_K$ de $K$ sans facteur carré, étranger à $\ell$ et stable par $G=\Gal(K/\k)$, pour chaque classe $[\wi\d_K]$ d'ordre $\ell$-primaire dans  $\,\wCl_K^{\,\m_{\si{K}}}$, il existe une infinité de $\ell$-extensions abéliennes $\infty$-décomposées $F/\k$ telles que la classe $[\wi\d_K]$ se principalise dans le groupe de classes de rayons $\,\wCl_L^{\,\m_{\si{L}}}$ du compositum $L=FK$.
\end{Th}

La première étape consiste à se ramener au cas où $[\wi\d_K]$ est une puissance $\ell^{m_{\si{K}}}$-ième dans $\,\wCl_K^{\,\m_{\si{K}}}$:

\begin{Prop}\label{Puissance}
Sous la conjecture de Gross-Kuz'min, quitte à grossir $K$ par composition avec une $\ell$-extension abélienne $\infty$-décomposée de $\k$, on peut supposer $[\wi\d_K] \in (\wCl_K^{\,\m_{\si{K}}})^{\ell^m}$ avec $m$ arbitaire.
\end{Prop}

\Preuve Notons $K^c$ la $\Zl$-extension cyclotomique de $K$; puis $K_{\scale{0.8}{\infty}}=K[\zeta_{\ell^{^{\scale{0.6}{\infty}}}}]=K^c[\zeta_{2\ell}]$ l'extension engendrée par les racines d'ordre $\ell$-primaire de l'unité; et considérons la pro-$\ell$-extension $\wH^{\,\m_{\si{K}}}_K[\zeta_{2\ell}]$.

\begin{itemize}
\item Si $\ell$ est impair et $K$ vérifie la conjecture de Gross-Kuz'min pour $\ell$, le groupe de Galois $\Gal(\wH^{\,\m_{\si{K}}}_K[\zeta_{2\ell}]/K$ s'identifie au produit direct du groupe cyclique $\Delta=\Gal(\wH^{\,\m_{\si{K}}}_K[\zeta_{2\ell}]/\wH^{\,\m_{\si{K}}}_K)$, du $\ell$-groupe fini $\,\wCl_K^{\,\m_{\si{K}}} \simeq \Gal(\wH^{\,\m_{\si{K}}}_K/K^c)$ et du groupe procyclique $\Gamma=\Gal(K^c/K)\simeq\Zl$.
\item Si $\ell$ vaut $2$, toujours sous la conjecture de Gross-Kuz'min, il faut distinguer:
\begin{itemize}
\item Si $\wH^{\,\m_{\si{K}}}_K$ ne contient pas $i=\zeta_4$, la même décomposition vaut encore avec  $\Delta\simeq \ZZ/2\ZZ$ .
\item Si $\wH^{\,\m_{\si{K}}}_K$  contient $i$, l'extension $K[i]/K$ est localement cyclotomique; et en remplaçant $K$ par $K[i)=K\QQ[i]$, qui vérifie les mêmes hypothèses, on est ramené au cas précédent.
\end{itemize}
\end{itemize}
Soit alors $\wK_m/K$ la sous-extension d'exposant $\ell^m$ de  $\wH^{\,\m_{\si{K}}}_K[\zeta_{2\ell}]/K$. Dans tous les cas, le théorème de \v Cebotarev  appliqué dans la clôture galoisienne de $\wK_m$ sur $\QQ$ nous assure l'existence d'un premier $\q_K^\ph \nmid \ell\m^\ph_K$ de $K$  au-dessus d'un premier $q$ de $\NN$ complètement décomposé dans $K[\zeta_{2\ell^h}]/\QQ $ de même image que $\d_K^\ph$ dans $\wCl_K^{\,\m_{\si{K}}}/\wCl_K^{\,\m_{\si{K}}}{}^{\ell^m}$ (en notations multiplicatives).

Maintenant, le sous-corps réel du corps cyclotomique $\QQ[\zeta_{2\ell^m}]$ contient un unique sous-corps $F_\circ$ cyclique de degré $\ell^m$ et totalement ramifié en $q$. De plus, comme $q$ est pris complètement décomposé dans $K/\QQ$, la place $\q_K^\ph$ est totalement ramifiée dans l'extension composée $KF_\circ/K$; de sorte que l'étendue $[\wi\d_{K\!F_\circ}]$ de $[\wi\d_K^\ph]$  est bien une puissance $\ell^m$-ième dans $\,\wCl_{K\!F_\circ}^{\,\m_{\si{K\!F_\circ}}}$.\medskip

\PreuveTh Revenons aux notations additives. D'après la Proposition \ref{Puissance}, nous pouvons supposer $[\wi\d_K] \in \ell^{m_{\si{K}}}\,\wCl_K^{\,\m_{\si{K}}}$ sans restreindre la généralité, en notant $\ell^{m_{\si{K}}}$ l'ordre du $\ell$-groupe $\mu^\ph_K$. La Proposition \ref{Construction} nous assure alors pour $n$ arbitrairement grand l'existence d'une infinité de $\ell$-extensions cycliques $\infty$-décomposées $F$ de $\k$ logarithmiquement ramifiées en un unique premier $\p$ avec pour indice $\wi e_\p(F/\k)=\ell^n$, qui satisfont les conditions $(i)$, $(ii)$ et $(iii)$ de la section précédente.

En particulier la classe $[\wi\d^\ph_K]$ de $\d_K$ dans $\,\wCl_K^{\,\m_{\si{K}}}$ est représentée modulo ${\ell^n}\,\wCl_K^{\,\m_{\si{K}}}$ par l'un des $[K:\k]$ idéaux $\p^\ph_K$ au-dessus de $\p$, lequel se ramifie logarithmiquement dans l'extension composée $L/K=KF/K$ avec pour indice $\ell^n$, de sorte qu'il vient: 
$[\,\wi\d_K^\ph]=[\,\wi\p_K^\ph]+\ell^n  [\,\wi\r_K^\ph]$ et $\wi\p^\ph_K=\ell^n\,\wi\q_L^\ph$ pour un diviseur logarithmique $\wi \q^\ph_L$  invariant par $C=\Gal(L/K)$.
D'après la Proposition \ref{CdA} l'étendue $[\wi\d^\ph_L]= \ell^n\, [\wi\q^\ph_L + \wi\r^\ph_K]$ de $[\d^\ph_K]$ à $L$ est ainsi la classe principale dès qu'on a: $n\ge \wi c_{K,\m}$.

\begin{Sco}
La conclusion du Théorème vaut encore lorsque $K/\k$ n'est pas supposée galoisienne.
\end{Sco}

\Preuve Partons d'une classe d'ordre $\ell$-primaire $[\wi\d^\ph_K]$; introduisons la clôture galoisienne $\bar K/\k$ de $K/\k$ et notons $\ell^h$ la $\ell$-partie du degré de $\bar K/K$. D'après la Proposition \ref{Puissance}, quitte à grossir $K$ par composition avec une $\ell$-extension abélienne de $\k$ (ce qui n'augmente pas $h$), nous pouvons supposer que $[\wi\d^\ph_K]$ est une puissance $\ell^h$-ième dans $\,\wCl_K^{\,\m_{\si{K}}}$, donc la norme dans $\bar K/K$ d'une classe $[\wi\d^\ph_{\bar K}]$ de $\,\Cl_{\bar K}^{\,\m_{\si{\bar K}}}$. Donnons-nous maintenant un corps principalisant $F$ pour $[\d^\ph_{\bar K}]$, i.e. une $\ell$-extension abélienne $\infty$-décomposée $F/\k$ telle que l'étendue $[\wi\d^\ph_{F\bar K}]$ de $[\wi\d^\ph_{\bar K}]$ soit la classe triviale. Alors l'étendue à $FK$ de sa norme $N_{{\bar K/(\bar K \cap FK)}}([\wi\d^\ph_{\bar K}])$ est la classe triviale de $\,\wCl_{FK}^{\,\m_{\si{F\!K}}}$. Or on a:

\centerline{$[\wi\d^\ph_{ K}]=N_{\bar K/K}([\wi\d^\ph_{\bar K}])=N_{(\bar K\cap FK)/K}\big( N_{\bar K/(\bar K \cap FK)}(([\wi\d^\ph_{\bar K}])\big)$;}

\noindent et $[\wi\d^\ph_K]$ se principalise dans $\,\wCl_{FK}^{\,\m_{\si{F\!K}}}$.	

\newpage
\section{Retour sur la conjecture de Gross-Kuz'min}

L'objet de cette section est de mieux cerner le rôle exact de la conjecture de Gross-Kuz'min dans les diverses étapes de la preuve du résultat principal sur la principalisation. À l'analyse, dans celles-ci, la conjecture intervient essentiellement pour trois points-clés:\smallskip

\noindent{\em (i) Pour assurer la finitude des $\ell$-groupes de classes logarithmiques de degré nul}\medskip

Bien entendu, la capitulation ne peut concerner que des classes d'ordre fini. Si donc l'on veut se dispenser de la conjecture de Gross-Kuz'min, il convient de remplacer le pro-$\ell$-groupe des classes logarithmiques de rayons $\,\wCl_K^{\,\m_{\si{K}}}$ par son sous-groupe de torsion, disons $\,\wT_K^{\,\m_{\si{K}}}$, qui est toujours fini.

Dans ce cas, la sous-extension de $\wH_K^{\,\m_{\si{K}}}$ fixée par  $\,\wT_K^{\,\m_{\si{K}}}$ n'est plus, a priori, la $\Zl$-extension cyclotomique $K^c$ de $K$, mais simplement le compositum des $\Zl$-extensions localement cyclotomiques. Son groupe de Galois sur $K^c$ est un $\Zl$-module libre dont la dimension mesure précisément le défaut de la conjecture de Gross-Kuz'min dans le corps $K$ (pour le premier $\ell$ fixé).\medskip

\noindent{\em (ii) Pour construire le sous-module $\,\wE^{\,\varepsilon}_K$}\smallskip

L'existence d'un sous-module monogène  $\,\wE^{\,\varepsilon}_K$ isomorphe à l'idéal d'augmentation sous l'hypothèse de décomposition d'une place à l'infini repose sur l'expression du caractère des unités logarithmiques sous la conjecture de Gross-Kuz'min dans $K$. Mais, ici encore, il n'est nullement nécessaire de s'appuyer sur la conjecture: il est bien connu, en effet, que le groupe $\,\wE_K^\ph$ contient un sous-groupe canonique, noté $\,\wE_K^\nu$ dans \cite{J55} \S7, dit groupe des normes universelles de Kuz'min, qui possède le rang et donc le caractère attendus. De ce fait, le sous-groupe d'indice fini $\,\wE_K^\nu \cap \,\wE_K^{\,\m_{\si{K}}}$ contient bien un sous-module de type $\,\wE^{\,\varepsilon}_K$ indépendamment de la conjecture de Gross-Kuz'min.\medskip

\noindent{\em (ii) Pour évaluer le quotient de Herbrand $q(C,\,\wE_K^{\,\m_{\si{K}}})$}\smallskip

Ce point est le plus sensible: les identités des classes logarithmiques d'ambiges dans l'extension cyclique $L/K$ de groupe de Galois $C$ font apparaître le quotient 
$\Pl_L^{\,\m_{\si{L}}}{}^C/\Pl^{\,\m_{\si{K}}}_K \simeq  H^1(C,\wE^{\,\m_{\si{L}}}_L)$ dont l'ordre coïncide, sous la conjecture de Gross-Kuz'min {\em dans $L$}, avec celui du groupe $H^2(C,\wE^{\,\m_{\si{L}}}_L)$, le quotient de Herbrand $q(C, \wE_L^{\,\m_{\si{L}}})=q(C,\wE_L)$ du groupe d'unités logarithmiques $\,\wE^{\,\m_{\si{L}}}_L$ valant alors 1.

Or, ce quotient ne dépend que du caractère du groupe $\,\wE_L$. Si le groupe $C$ est cyclique d'ordre, disons, $\ell^h$, la décomposition semi-simple de l'algèbre de groupe\smallskip

\centerline{$\Ql[C] \simeq \Ql[X]/(X^{\ell^h}-1) \simeq \Ql \oplus \Ql[\zeta_\ell] \oplus \cdots \oplus \Ql[\zeta_{\ell^h}]$}\smallskip

\noindent conduit à la décomposition $\ell$-adique irréductible du caractère régulier\smallskip 

\centerline{$\chi_C^{\textrm{rég}}= \chi_0^\ph + \chi_1^\ph + \cdots + \chi_h^\ph =  \chi_0^\ph + \chi_C^{\textrm{aug}}$,}\smallskip

\noindent où $\chi_0^\ph$ désigne le caractère unité et $\chi_C^{\textrm{aug}}$ le caractère d'augmentation. Maintenant, pour tout $\Zl[C]$-module projectif $M$ de caractère $\chi_M^\ph=\sum n_i\chi_i$, un calcul élémentaire montre que l'ordre $h^2(M)$ du groupe $H^2(C,M)$ ne dépend que de $n_0$, tandis que l'ordre $h^1(M)$ du groupe $H^1(C,M)$ ne dépend que des $n_i$ pour $i>0$. Écrivons donc $\,\wE'_L=\,\wE_L^\ph/\wE_L^\nu$ le quotient du groupe des unités logarithmiques par le sous-groupe des normes universelles de Kuz'min introduit plus haut.

Sous la conjecture de Gross-Kuz'min {\em dans $K$}, le groupe  $\,\wE'_K=\,\wE_K^\ph/\wE_K^\nu$ est fini, de sorte que la composante unité du caractère de $\,\wE'_L$ est alors triviale. En particulier, son quotient de Herbrand $q(C,\,\wE'_L) =h^1(\,\wE'_L)/h^2(\,\wE'_L)$ est une puissance positive de $\ell$. Et il vient donc:\smallskip

\centerline{$h^1(\,\wE^{\,\m_{\si{L}}}_L)=h^2(\,\wE^{\,\m_{\si{L}}}_L)q(\,\wE^{\,\m_{\si{L}}}_L)=h^2(\,\wE^{\,\m_{\si{L}}}_L)q(\,\wE_L^\ph)=h^2(\,\wE^{\,\m_{\si{L}}}_L)q(\,\wE'_L) \ge h^2(\,\wE^{\,\m_{\si{L}}}_L)$,}\smallskip

\noindent en vertu de l'identité $q(C,\,\wE_L^{\,\nu})=1$ valable inconditionnellement.\medskip

En fin de compte, dès lors que le seul corps $K$ vérifie la conjecture de Gross-Kuz'min pour $\ell$, les calculs de la Proposition \ref{CLA} fournissent la {\em majoration} du nombre de classes d'ambiges:
$$
(\wDl_L^{\,\m_{\si{L}}}{}^C:\Pl_L^{\,\m_{\si{L}}}{}^C) \;\le \; |\,\wCl_K^{\,\m_{\si{K}}}|\;\frac{\prod_{\p_{\si{\!\infty}}} d_{\p_{\si{\!\infty}}}\!(L/K)\,\prod_{\p_{\si{\circ}}\nmid\m} \tilde e_{\p_{\si{\circ}}}(L/K)}{\big(\deg_L\wDl_L^{\,\m_{\si{L}}}{}^C:\deg_L\wDl^{\,\m_{\si{K}}}_K\big)\;\big( \wE^{\,\m_{\si{K}}}_K:N_{L/K}(\wE^{\,\m_{\si{L}}}_L) \big)}
$$ 
Et la conclusion de la Proposition \ref{CdA} reste valable; ce qui achève la démonstration du Théorème \ref{ThP}.

\newpage
\section{Conséquences arithmétiques}

Dans toute cette section, nous supposons fixé le nombre premier $\ell$.

\begin{Def}
Soient $\k$ un corps de nombres et $T=T_\k$ un ensemble fini de premiers de $\k$ ne divisant pas $\ell$. Pour toute extension algébrique $N$ de $\k$ (non nécessairement de degré fini sur $\k$), convenons de noter $T_N$ l'ensemble (éventuellement infini) des places de $N$ au-dessus de $T$ et, pour $[K:\k]$ fini, par $\m_K$ le diviseur  $\m_K^\ph=\prod_{\q_K^\ph\in T_K^\ph}\q_K^\ph$.
Nous disons que $N$ est $T_N^{\si{\,(\ell)}}$-logarithmiquement principal lorsque son groupe des classes $T_N^{\si{\,(\ell)}}$-infinitésimales 
(défini comme limite inductive des $\ell$-groupes logarithmiques de classes de rayons $\,\wCl_K^{\,\m_{\si{K}}}$ associés aux sous-extensions $K$ de degré fini) est trivial; ce qui s'écrit:  
 
 \centerline{$\wCl_K^{T_N^{\si{\,(\ell)}}}=\varinjlim \,\wCl_K^{\,\m_{\si{K}}}=1$.}\smallskip
 
 Lorsque $T$ est vide, nous disons tout simplement que $N$ est logarithmiquement principal.
\end{Def}

Cela posé, appliquons d'abord le Théorème principal avec $\k=\QQ$ logarithmiquement principal.

\begin{Sco}\label{Sco}
Pour tout corps de nombres totalement réel $K$ vérifiant la conjecture de Gross-Kuz'min en $\ell$ (par exemple un corps abélien réel)  et tout ensemble fini $T=T_\QQ$ de nombres premiers $p \ne \ell$, il existe une infinité d'extensions abéliennes réelles $F/\QQ$ telles que les $\ell$-classes logarithmiques de rayons modulo $\m_K^\ph=\prod_{\q_K^\ph\in T_K^\ph}\q_K^\ph$ de $\,\wCl_K^{\,\m_{\si{K}}}$ capitulent dans  $\,\wCl_{F\!K}^{\,\m_{\si{F\!K}}}$.\par
Autrement dit, le groupe  $\,\wCl_K^{\,\m_{\si{K}}}$  capitule dans le sur-corps $K[\cos(2\pi/n)]$ pour une infinité de $n$.
\end{Sco}

Passant à la limite inductive dans le cas particulier $T=\emptyset$, nous obtenons:

\begin{Cor}
Le sous-corps réel maximal $\QQ^\ab_+=\bigcup_{n>0}\QQ[\cos(2\pi/n)]$ du corps cyclotomique $\QQ^\ab=\bigcup_{n>0}\QQ[\zeta_n]$ est logarithmiquement principal.\par
Plus généralement les extensions algébriques totalement réelles $N$ de $\QQ^\ab_+$ qui vérifient la conjecture de Gross-Kuz'min en $\ell$ sont logarithmiquement principales.
\end{Cor}

Regardons maintenant le cas relatif en distinguant suivant la signature du corps de base $\k$:

\begin{Sco}
Soit $\k$ un corps de nombres totalement réel et $T_\k$ un ensemble fini de places de $\k$ ne divisant pas $\ell$.
Alors, pour tout corps de nombres $K$ qui contient $\k$, satisfait la conjecture de Gross-Kuz'min en $\ell$ et possède au moins une place réelle, il existe une infinité d'extensions abéliennes totalement réelles $F/\k$ telles que les classes logarithmiques de rayons modulo $\m_K^\ph=\prod_{\q_K^\ph\in T_K^\ph}\q_K^\ph$ contenues dans le sous-groupe relatif $\,\wCl_{K/\k}^{\,\m_{\si{K}}}$ capitulent dans $\,\wCl_{F\!K}^{\,\m_{\si{F\!K}}}$.
\end{Sco}

\begin{Cor}
Soit $\k$ un corps de nombres totalement réel,  $\k^\ab_+$ sa plus grande extension abélienne totalement réelle et $T_\k$ un ensemble fini de places de $\k$ qui ne divisent pas $\ell$.\par
Si $\k$ est $T_\k^{\si{\,(\ell)}}$-logarithmiquement principal, toute extension algébrique $K$ de $k^\ab_+$ qui satisfait la conjecture de Gross-Kuz'min et dont la clôture galoisienne possède au moins un plongement réel est encore $T_K^{\si{\,(\ell)}}$-logarithmiquement principale.
\end{Cor}

\Preuve Pour tout $\alpha$ dans $K$, le groupe de classes logarithmiques $\,\wCl_{\k[\alpha]}^{\,\m_{\si{\k[\alpha]}}}$ du corps $\k[\alpha]$ capitule dans le sous-corps $\k^\ab_+[\alpha]$ de $K$. Et $K$ est donc bien $T_K^{\si{\,(\ell)}}$-logarithmiquement principal, comme annoncé.
\begin{Sco}
Soit $\k$ un corps de nombres qui possède au moins une place complexe (par exemple un corps quadratique imaginaire) 
et $T_\k$ un ensemble fini de places de $\k$ ne divisant pas $\ell$.
Alors, pour tout corps de nombres $K$ qui contient $\k$, il existe une infinité d'extensions abéliennes $F/\k$ complètement décomposées en toutes les places à l'infini et telles que les classes logarithmiques de rayons modulo $\m_K^\ph=\prod_{\q_K^\ph\in T_K^\ph}\q_K^\ph$ contenues dans $\,\wCl_{K/\k}^{\,\m_{\si{K}}}$ capitulent dans $\,\wCl_{F\!K}^{\m_{\si{F\!K}}}$.
\end{Sco}

\begin{Cor}
Soient $\k$ un corps de nombres qui possède au moins une place complexe (par exemple un corps quadratique imaginaire), $\k^\ab_+$ la plus grande extension abélienne de $\k$ complètement décomposées en toutes les places à l'infini, et $T_\k$ un ensemble fini de places de $\k$.\par

Si $\k$ est $T_\k^{\si{\,(\ell)}}$-logarithmiquement principal, toute extension algébrique $K$ de $\k^\ab_+$ qui satisfait la conjecture de Gross-Kuz'min est encore $T_K^{\si{\,(\ell)}}$-logarithmiquement principale.
\end{Cor}

\newpage
\noindent{\sc  Appendice}: {\bf Parallèle avec la conjecture de Grenberg}\medskip

Il peut être intéressant de retranscrire le Théorème principal de \cite{J60} sous la forme:

\begin{ThA} 
Soit $K$ un corps de nombres totalement réel et $\ell$ un nombre premier donné. Sous la conjecture de Leopoldt dans $K^c$  (par exemple pour $K$ abélien), la conjecture de Greenberg pour $K$ en $\ell$ revient à postuler que  le $\ell$-groupe  $\,\wCl_K$ des classes logarithmiques du corps $K$ capitule dans le corps totalement réel $K[\cos(2\pi/\ell^m]$ pour tout $m$ assez grand.
\end{ThA}

Cette formulation est alors à mettre en parallèle avec le Scolie \ref{Sco} plus haut, qui implique:

\begin{ThB}
Si $K$ est un corps de nombres totalement réel qui vérifie la conjecture de Gross-Kuz'min pour un premier $\ell$ (par exemple si $K$ vérifie la conjecture de Leopoldt en $\ell$, notamment si $K$ est abélien), le $\ell$-groupe  $\,\wCl_K$ capitule dans le corps totalement réel $K[\cos(2\pi/n]$ pour une infinité de $n$ qui ne sont pas divisibles par $\ell$ (en fait qui peuvent être pris étrangers à tout entier donné à l'avance).
\end{ThB}

\addcontentsline{toc}{section}{Bibliographie}
\def\refname{\normalsize{\sc  Références}}
\bigskip
{\small
\noindent{\sc Adresse:}
Univ. Bordeaux \& CNRS,
Institut de Mathématiques de Bordeaux, UMR 5251,
351 Cours de la Libération,
F-33405 Talence cedex

\noindent{\sc Courriel:}
 \tt jean-francois.jaulent@math.u-bordeaux.fr
}

\end{document}